\documentstyle[11pt,twoside]{article}
\newfont{\Bbb}{msbm10 scaled\magstephalf}
\begin{document}
\title{ Weighted Composition Operators between different Bloch-type Spaces in Polydisk }
\author{{\small\sc Zehua Zhou\hspace*{6mm}Renyu Chen}\\
Department of Mathematics, Tianjin University\\
Tianjin 300072, P. R. China\\ E-mail: zehuazhou2003@yahoo.com.cn}
\date{}
\maketitle \footnote{Zehua Zhou, Corresponding author. Supported
in part by the National Natural Science Foundation of China (Grand
No.10371091), and LiuHui Center for Applied Mathematics, Nankai
University and Tianjin University}

\begin{abstract}Let $\phi(z)=(\phi_1(z), \ldots,\phi_n(z))$ be a holomorphic
self-map of $U^n$ and $\psi(z)$ a holomorphic function on $U^n,$
where $U^n$ is the unit polydisk of $\mbox{\Bbb C}^n.$ Let $p\geq
0,$ $q\geq 0$, this paper gives some necessary and sufficient
conditions for the weighted composition operator $W_{\psi,\phi}$
induced by $\psi$ and $\phi$ to be bounded and compact between
$p$-Bloch space ${\cal B}^p(U^n)$ and $q$-Bloch space ${\cal
B}^q(U^n).$
\end{abstract}

\hspace{4mm}{\bf Keywords} \hspace{2mm}Bloch-type space,
 Weighted composition operator, Composition operator.

\hspace{4mm}2000\hspace{1mm}Mathematics Subject Classification
\hspace{2mm}47B38, 47B33, 32A37, 32A30

\section{Introduction}
\newtheorem{Definition}{Definition}
\newtheorem{Lemma}{Lemma}
\newtheorem{Theorem}{Theorem}
\newtheorem{Corollary}{Corollary}

Let $\Omega$ be a bounded homogeneous domain in $\mbox{\Bbb
C}^{n}.$ The class of all holomorphic functions with domain
$\Omega$ will be denoted by $H(\Omega).$ Let $\phi $ be a
holomorphic self-map of $\Omega,$ the composition operator
$C_{\phi}$ induced by $\phi$ is defined by
$$(C_{\phi}f)(z)=f(\phi(z)),$$ for $z$ in $\Omega$ and $f\in
H(\Omega)$. If, in addition, $\psi$ is a holomorphic function
defined on $\Omega,$ the multiplication operator induced by $\psi$
is defined by $$M_{\psi}f(z)=\psi(z)f(z),$$ and weighted
composition operators $W_{\psi,\phi}$ induced by $\psi$ and $\phi$
is defined by
$$\left(W_{\psi,\phi}f\right)(z)=\psi(z)f(\phi(z))$$ for $z$ in
$\Omega$ and $f\in H(\Omega).$  If let $\psi\equiv 1,$ then
$W_{\psi,\phi}=C_{\phi},$  if let $\phi=Id$, then
$W_{\psi,\phi}=M_{\psi}.$ So we can regard weighted composition
operator as a generalization of a multiplication operator and a
composition operator.

Let $K(z,z)$ be the Bergman kernel function of $\Omega$, the
Bergman metric $H_{z}(u,u)$ in $\Omega$ is defined by
$$H_{z}(u,u)=\displaystyle\frac{1}{2}
\sum\limits^{n}_{j,k=1} \displaystyle\frac{\partial^{2}\log
K(z,z)}{\partial z_{j}
\partial {\overline {z}}_{k}}u_{j}{\overline u}_{k},$$
where $z\in\Omega$ and $u=(u_{1},\ldots,u_{n})\in \mbox{\Bbb
C}^{n}.$

Following Timoney [1], we say that $f\in H(\Omega)$ is in the
Bloch space ${\cal B}(\Omega),$ if
$$\|f\|_{{\cal B}(\Omega)}=\sup\limits_{z\in \Omega}Q_{f}(z)<\infty,$$
where $$Q_{f}(z)=\sup\left\{\displaystyle\frac {|\bigtriangledown
f(z)u|}{H^{\frac{1}{2}}_{z}(u,u)}: u\in \mbox{\Bbb
C}^{n}-\{0\}\right\},$$ and $\bigtriangledown f(z)
=\left(\frac{\partial f(z)}{\partial z_{1}}, \ldots,
\frac{\partial f(z)}{\partial z_{n}} \right), \bigtriangledown
f(z)u =\sum\limits^{n}_{l=1}\frac{\partial f(z)}{\partial
z_{l}}u_{l}.$

The little Bloch space ${\cal B}_0(\Omega)$ is the closure in the
Banach space ${\cal B}(\Omega)$ of the polynomial functions.

Let $\partial\Omega$ denote the boundary of $\Omega$. Following
Timoney [2], for $\Omega=B_n$ the unit ball of $\mbox{\Bbb C}^n$,
${\cal B}_0(B_n)=\left\{f\in{\cal B}(B_n): Q_f(z)\to 0,
\mbox{as}\hspace*{2mm}z\to\partial B_n \right\};$ for $\Omega=\cal
D$ the bounded symmetric domain other than $B_n$, $\left\{
f\in{\cal B}({\cal D}): Q_f(z)\to 0, \mbox{as}\hspace*{2mm}
z\to\partial{\cal D}\right\}$ is the set of constant functions on
$\cal D.$ So if $\cal D$ is a bounded symmetric domain other than
the ball, we denote the ${\cal B}_{0*}({\cal D})= \left\{f\in{\cal
B}({\cal D}): Q_f(z)\to 0, \mbox{as}\hspace*{2mm}
z\to\partial^*{\cal D}\right\}$ and also call it little Bloch
space, here $\partial^*{\cal D}$ means the distinguished boundary
of $\cal D$. The unit ball is the only bounded symmetric domain
$\cal D$ with the property that $\partial^*{\cal D}=\partial{\cal
D}.$

Let $U^n$ be the unit polydisk of $\mbox{\Bbb C}^n.$ Timony [1]
shows that $f\in(U^n)$ if and only if
$$\|f\|_=|f(0)|+\sup\limits_{z\in U^n}\sum\limits^n_{k=1}
\left|\displaystyle\frac{\partial f} {\partial
z_k}(z)\right|\left(1- |z_k|^2\right)<+\infty.$$ This equality was
the starting point for introducing the $p$-Bloch spaces.

Let $p\geq 0,$ a function $f\in H(U^n)$ is said to belong to the
$p$-Bloch space ${\cal B} ^p(U^n)$ if
$$\|f\|_p=|f(0)|+\sup\limits_{z\in U^n}\sum\limits^n_{k=1}
\left|\frac{\partial f} {\partial z_k}(z)\right|\left(1-
|z_k|^2\right)^p<+\infty.$$ It is easy to show that ${\cal
B}^p(U^n)$ is a Banach space with the norm $\|\cdot\|_p.$

It is easy to see that if
$$\lim_{z\to\partial U^n}\sum\limits^n_{k=1}\left|\frac{\partial f}{\partial
z_k}(z)\right|(1- |z_k|^2)^p=0,$$ then $f$ must be a constant. So,
there is no sense to introduce the corresponding little $p$-Bolch
space in this way. We will say that the little $p$-Bolch space
${\cal B}_0^p(U^n)$ is the closure of the polynomials in the
$p$-Bolch space.

In the recent years, there have been many papers focused on
studying the composition operators in function spaces (say, for
1-dimensional case see [3-8], for n-dimensional case see [9-12]).
More recently, S.Ohno, K. Stroethoff and R.H.Zhao in [8] discuss
the weighted composition operators between-type spaces for
1-dimensional case.

In this paper, we discuss the boundedness and compactness of the
weighted composition operators between $p$-Bloch space and
$q$-Bloch space in polydisk, some new methods and techniques have
been used because of the difference between topology boundary
$\partial U^n$ and distinguished boundary $\partial^* U^n$ of
$U^n$, where $n>1$, especially in the proof of Theorems 2 and 3
(need to be discussed according to the properties of the
boundary). The results in this paper will extend corresponding
results on the Bloch spaces (see [4-11]).

Our main results are the following:
\begin{Theorem} Let $\phi=(\phi_1, \ldots, \phi_n)$ be a holomorphic self-map
of $U^n$ and $\psi(z)$ a holomorphic function of
$U^n,$ $p\geq 0,$ $q\geq 0.$ \\
(i) If $p=1$, then $W_{\psi,\phi}: {\cal B}^p(U^n)\longrightarrow
{\cal B}^q(U^n) $ is bounded if and only if
\begin{equation}\sum\limits^n_{k,l=1}\left|\displaystyle\frac{\partial \psi}
{\partial z_k}(z)\right|\left(1-|z_k|^2\right)^{q}\ln
\displaystyle\frac{4}{1-|\phi_l(z)|^2}=O(1) \hspace*{6mm}(z\in U^n
)\label{1}\end{equation} and
\begin{equation}|\psi(z)|\sum\limits^n_{k,l=1}\left|\displaystyle\frac{\partial \phi_{l}}
{\partial
z_k}(z)\right|\displaystyle\frac{\left(1-|z_k|^2\right)^{q}}{1-|\phi_l(z)|^2}=O(1)
\hspace*{6mm}( z\in U^n ).\label{2}\end{equation} (ii) If $0\leq
p<1$, then $W_{\psi,\phi}: {\cal B}^p(U^n)\longrightarrow {\cal
B}^q(U^n) $ is bounded if and only if
\begin{equation}\psi\in {\cal B}^q(U^n)\label{3}\end{equation} and
\begin{equation}|\psi(z)|\sum\limits^n_{k,l=1}\left|\displaystyle\frac{\partial
\phi_{l}} {\partial
z_k}(z)\right|\displaystyle\frac{\left(1-|z_k|^2\right)^{q}}{\left(1-|\phi_l(z)|^2\right)^{p}}=O(1)
\hspace*{6mm}( z\in U^n ).\label{4}\end{equation} (iii) If $p>1$,
then $W_{\psi,\phi}: {\cal B}^p(U^n)\longrightarrow {\cal
B}^q(U^n) $ is bounded if and only if
\begin{equation}\sum\limits^n_{k,l=1}\left|\displaystyle\frac{\partial \psi}
{\partial
z_k}(z)\right|\displaystyle\frac{\left(1-|z_k|^2\right)^{q}}{\left(1-|\phi_l(z)|^2\right)^{p-1}}
=O(1) \hspace*{6mm}( z\in U^n)\label{5}\end{equation} and
\begin{equation}|\psi(z)|\sum\limits^n_{k,l=1}\left|\displaystyle\frac{\partial
\phi_{l}} {\partial
z_k}(z)\right|\displaystyle\frac{\left(1-|z_k|^2\right)^{q}}{\left(1-|\phi_l(z)|^2\right)^{p}}=O(1)
\hspace*{6mm}( z\in U^n ).\label{6}\end{equation}
\end{Theorem}

\begin{Theorem} Let $\phi=(\phi_1, \ldots, \phi_n)$ be a holomorphic self-map of $U^n$
and $\psi(z)$ a holomorphic function of $U^n$,$p\geq 0,
q\geq 0.$\\
(i) If $p=1$, then $W_{\psi,\phi}: {\cal B}^p(U^n)\longrightarrow
{\cal B}^q(U^n) $ is compact if and only if $W_{\psi,\phi}: {\cal
B}^p(U^n)\longrightarrow {\cal B}^q(U^n)$ is bounded and
\begin{equation}\sum\limits^n_{k,l=1}\left|\displaystyle\frac{\partial \psi}
{\partial z_k}(z)\right|\left(1-|z_k|^2\right)^{q}\ln
\displaystyle\frac{4}{1-|\phi_l(z)|^2}=o(1)
\hspace*{6mm}(\phi(z)\to\partial U^n),\label{7}\end{equation}
and\begin{equation}|\psi(z)|\sum\limits^n_{k,l=1}\left|\displaystyle\frac{\partial
\phi_{l}} {\partial
z_k}(z)\right|\displaystyle\frac{\left(1-|z_k|^2\right)^{q}}{1-|\phi_l(z)|^2}
=o(1)\hspace*{6mm}(\phi(z)\to\partial U^n).\label{8}\end{equation}
(ii) If $p>1$, then $W_{\psi,\phi}: {\cal B}^p(U^n)\longrightarrow
{\cal B}^q(U^n) $ is compact if and only if $W_{\psi,\phi}: {\cal
B}^p(U^n)\longrightarrow {\cal B}^q(U^n)$ is bounded and
\begin{equation}\sum\limits^n_{k,l=1}\left|\displaystyle\frac{\partial \psi}
{\partial
z_k}(z)\right|\displaystyle\frac{\left(1-|z_k|^2\right)^{q}}{\left(1-|\phi_l(z)|^2\right)^{p-1}}
=o(1)\hspace*{6mm}(\phi(z)\to\partial U^n),\label{9}\end{equation}
and
\begin{equation}|\psi(z)|\sum\limits^n_{k,l=1}\left|\displaystyle\frac{\partial
\phi_{l}} {\partial
z_k}(z)\right|\displaystyle\frac{\left(1-|z_k|^2\right)^{q}}{\left(1-|\phi_l(z)|^2\right)^{p}}
=o(1)\hspace*{6mm}(\phi(z)\to\partial
U^n).\label{10}\end{equation}
\end{Theorem}

{\bf Remark 1} It is easy to show that if (1) or (5) or (7) or (9)
holds, then $\psi\in{\cal B}^q(U^n)$ and
\begin{equation}\sum\limits^n_{k=1}\left|\displaystyle\frac{\partial \psi(z)}
{\partial z_k}(z)\right|\left(1-|z_k|^2\right)^{q}
=o(1)\hspace*{6mm}( \phi(z)\to\partial
U^n).\label{11}\end{equation} In fact, (1) implies that
$$\sum\limits^n_{k}\left|\displaystyle\frac{\partial \psi}
{\partial
z_k}(z)\right|\left(1-|z_k|^2\right)^{q}\sum\limits^n_{l=1}\ln
\displaystyle\frac{4}{1-|\phi_l(z)|^2}\leq C$$ for all
$z\in\mbox{\Bbb C}^n.$ The same reason for others.

\begin{Theorem} Let $\phi=(\phi_1, \ldots, \phi_n)$ be a
holomorphic self-map of $U^n$ and $\psi(z)$ a holomorphic function
of $U^n$, $0\leq p<1,$ $q\geq 0.$ Then $W_{\psi,\phi}: {\cal
B}^p(U^n)\longrightarrow {\cal B}^q(U^n) $ is compact if and only
if $W_{\psi,\phi}: {\cal B}^p(U^n)\longrightarrow {\cal B}^q(U^n)
$ is bounded and
\begin{equation}\lim\limits_{|\phi_l(z)|\to 1}\sum\limits^n_{k=1}\left|\frac{\partial
\phi_{l}} {\partial
z_k}(z)\right|\frac{(1-|z_k|^2)^q}{(1-|\phi_l(z)|^2)^p}=0\label{12}\end{equation}
for each $l\in\{1,...,n\}$ with $|\phi_l(z)|\to 1.$\end{Theorem}

\begin{Corollary} Let $\phi=(\phi_1,...,\phi_n)$ be a holomorphic
self-map of $U^n$, $p,q\geq 0.$ Then $C_{\phi}:{\cal B} ^p(U^n)\to
{\cal B} ^q(U^n)$ is bounded if and only if there exists a
constant $C$ such that
$$\sum\limits^n_{k,l=1}\left|\frac{\partial \phi_{l}} {\partial
z_k}(z)\right|\frac{(1-|z_k|^2)^q} {(1-|\phi_l(z)|^2)^p}\leq C ,$$
for all $z\in U^n.$\end{Corollary}

{\bf Proof}\hspace*{4mm}Let $\psi(z)\equiv 1, z\in U^n$, then
$\displaystyle\frac{\partial \psi} {\partial z_k}(z)=0$ for all
$k\in\{1,2,\cdots,n\}.$  It is clear that condition (\ref{1}) in
Theorem 1, condition (\ref{7}) in Theorem 2 and condition
(\ref{9}) in Theorem 3 hold,  note that $W_{\psi,\phi}=C_{\phi}$,
the Corollary follows by combining the above Theorems.

Similar to Corollary 1, the following Corollary follows.

\begin{Corollary} Let $\phi=(\phi_1,...,\phi_n)$ be a holomorphic
self-map of $U^n.$ If $p\geq 1$ and $q\geq 0,$ then
$C_{\phi}:{\cal B} ^p(U^n)\to {\cal B} ^q(U^n)$ is compact if and
only if
$$\sum\limits^n_{k,l=1}\left|\frac{\partial \phi_{l}} {\partial
z_k}(z)\right|\frac{(1-|z_k|^2)^q} {(1-|\phi_l(z)|^2)^p}\leq C ,$$
for all $z\in U^n,$ and
$$\sum\limits^n_{k,l=1}\left|\frac{\partial \phi_{l}} {\partial
z_k}(z)\right|\frac{(1-|z_k|^2)^q}{(1-|\phi_l(z)|^2)^p} =o(1)\quad
( \phi(z)\to \partial U^n ).$$ If $0\leq p<1)$ and $q\geq 0,$ Then
$C_{\phi}: {\cal B} ^p(U^n)\to {\cal B} ^q(U^n)$ is compact if and
only if
$$\sum\limits^n_{k,l=1}\left|\frac{\partial \phi_{l}} {\partial
z_k}(z)\right|\frac{(1-|z_k|^2)^q} {(1-|\phi_l(z)|^2)^p}\leq C ,$$
for all $z\in U^n,$ and
$$\lim\limits_{|\phi_l(z)|\to
1}\sum\limits^n_{k=1}\left|\frac{\partial \phi_{l}} {\partial
z_k}(z)\right|\frac{(1-|z_k|^2)^q}{(1-|\phi_l(z)|^2)^p}=0$$for
each $l\in\{1,...,n\}$ with $\phi_l(z)|\to 1.$\end{Corollary}

{\bf Remark 2} If let $\phi=id: U^n\longrightarrow U^n,$ then we
can obtain the corresponding results about multiplication
operator$M_{\psi}:{\cal B} ^p(U^n)\to {\cal B}^q(U^n).$

Throughout the remainder of this paper $C$ will denote a positive
constant, the exact value of which will vary from one appearance
to the next.

\section{Some Lemmas}

\begin{Lemma} Let $f\in {\cal B}^p(U^n),$

(i) If $p=1,$ then $|f(z)|\leq
\left(\displaystyle\frac{1}{2}+\displaystyle\frac{1}{2n\ln
2}\right)\left(\sum\limits^n_{l=1}\ln\displaystyle\frac{4}{1-|z_l|^2}\right)\|f\|_{p};$

(ii) If $0\leq p<1,$ then $|f(z)|\leq
\left(1+\displaystyle\frac{1}{1-p}\right)\|f\|_{p};$

(iii) If $p>1$, then $|f(z)|\leq
\left(\displaystyle\frac{1}{n}+\displaystyle\frac{2^{p-1}}{p-1}\right)\left(\sum\limits^n_{l=1}
\displaystyle\frac{1}{(1-|z_l|^2)^{p-1}}\right)
\|f\|_{p}.$\end{Lemma}

{\bf Proof}\hspace{2mm} This Lemma can be proved by some integral
estimates (if necessary, the proof can be omitted).

By the definition of $\|.\|_{p}$,
$$|f(0)|\leq \|f\|_{p},\hspace*{4mm}\left|\displaystyle\frac{\partial f(z)
}{\partial z_l}\right|\leq
\displaystyle\frac{\|f\|_{p}}{(1-|z_l|^2)^p}
\hspace*{4mm}(l\in\{1,2,\cdots,n\})$$ and
\begin{eqnarray*}&&f(z)-f(0)=
\int^1_0 \displaystyle\frac{d f(tz)}{d
t}dt=\sum\limits^n_{l=1}\int^1_0 z_l\displaystyle\frac{\partial
f}{\partial\zeta_l}(tz)dt,\end{eqnarray*}
So\begin{eqnarray}&&|f(z)|\leq
|f(0)|+\sum\limits^n_{l=1}|z_l|\int^1_0\displaystyle\frac{\|f\|_p}{\left(1-t^2|z_l|^2\right)^p}dt
\nonumber\\
&&\leq
\|f\|_p+\|f\|_p\sum\limits^n_{l=1}\int^{|z_l|}_0\displaystyle\frac{1}{\left(1-t^2\right)^p}dt.\label{13}
\end{eqnarray}
If $p=1,$
\begin{equation}\int^{|z_l|}_0\displaystyle\frac{1}{\left(1-t^2\right)^p}dt=\displaystyle\frac{1}{2}
\ln\displaystyle\frac{1+|z_l|}{1-|z_l|}\leq
\displaystyle\frac{1}{2}
\ln\displaystyle\frac{4}{1-|z_l|^2}.\label{14}\end{equation} It is
clear that $\ln\displaystyle\frac{4}{1-|z_l|^2}>\ln 4=2\ln 2,$
so\begin{equation}1\leq \displaystyle\frac{1}{2\ln 2}
\ln\displaystyle\frac{4}{1-|z_l|^2}\leq
\displaystyle\frac{1}{2n\ln 2}
\sum\limits^n_{l=1}\ln\displaystyle\frac{4}{1-|z_l|^2}.\label{15}\end{equation}
Combining (\ref{13}),(\ref{14}) and (\ref{15}), we get
$$|f(z)|\leq \left(\displaystyle\frac{1}{2}+\displaystyle\frac{1}{2n\ln
2}\right)\left(\sum\limits^n_{l=1}\ln\displaystyle\frac{4}{1-|z_l|^2}\right)\|f\|_p.$$
If $p\neq 1,$
\begin{eqnarray}&&\int^{|z_l|}_0\displaystyle\frac{1}{\left(1-t^2\right)^p}dt
=\int^{|z_l|}_0\displaystyle\frac{1}{(1-t)^p}\cdot
\displaystyle\frac{1}{(1+t)^p}dt\nonumber\\
&&\leq
\int^{|z_l|}_0\displaystyle\frac{1}{(1-t)^p}dt=\displaystyle\frac{1-(1-|z_l|)^{-p+1}}{1-p}.\label{16}
\end{eqnarray}
If $0\leq p<1,$ (\ref{16}) gives that
$\int^{|z_l|}_0\displaystyle\frac{1}{\left(1-t^2\right)^p}dt\leq\displaystyle\frac{1}{1-p},$
it follows from (\ref{13}) that
$|f(z)|\leq\left(1+\displaystyle\frac{n}{1-p}\right)\|f\|_p.$

If $p>1,$ (\ref{16}) gives that
$$\int^{|z_l|}_0\displaystyle\frac{1}{\left(1-t^2\right)^p}dt\leq \displaystyle\frac{1-(1-|z_l|^{p-1})}{(p-1)(1-|z_l|)^{p-1}}
\leq\displaystyle\frac{2^{p-1}}{(p-1)(1-|z_l^2|)^{p-1}},$$ it
follows from (\ref{13}) that
\begin{eqnarray*}|f(z)|&\leq&\|f\|_p+\displaystyle\frac{2^{p-1}}{p-1}\left(\sum\limits^n_{l=1}
\displaystyle\frac{1}{(1-|z_l|^2)^{p-1}}\right)\|f\|_p\\
&\leq&
\left(\displaystyle\frac{1}{n}+\displaystyle\frac{2^{p-1}}{p-1}\right)\left(\sum\limits^n_{l=1}
\displaystyle\frac{1}{(1-|z_l|^2)^{p-1}}\right)
\|f\|_{p}.\end{eqnarray*}Now the Lemma is proved.

\begin{Lemma} Let $\phi$ be a holomorphic self-map of $U^n$ and $\psi(z)$ a holomorphic function
on $U^n,$ then $W_{\psi,\phi}:{\cal B} ^p(U^n)\to {\cal B}
^q(U^n)$ is compact if and only if for any bounded sequence
$\{f_j\}$ in ${\cal B} ^p(U^n)$ which converges to zero uniformly
on compact subsets of $U^n$, we have $\|W_{\psi,\phi}f_j\|_q\to
0,$ as $j\to\infty.$
\end{Lemma}

{\bf Proof}\hspace{4mm} Assume that $W_{\psi,\phi}$ is compact and
suppose $\{f_j\}$ is a sequence in ${\cal B} ^p(U^n)$ with
$\sup\limits_{j\in {\bf N}}\|f_j\|_p<\infty$ and $f_j\to 0$
uniformly on compact subsets of $U^n.$  By the compactness of
$W_{\psi,\phi}$ we have that $W_{\psi,\phi}(f_j)=\psi
f_j\circ\phi$ has a subsequence $\psi f_{j_m}\circ\phi$ which
converges in ${\cal B}^q,$ say, to $g.$ By Lemma 1 we have that
for any compact $K\subset U^n$ there is a positive constant $C_K$
independent of $f$ such that
$$|\psi(z)f_{j_m}(\phi(z))-g(z)|\leq C_K\|\psi f_{j_m}\circ \phi-g\|_q $$ for all $z\in
K$. This implies that $\psi(z)f_{j_m}(\phi(z))-g(z)\to 0$
uniformly on compacts of $U^n.$ Since $K$ is a compact subset of
$U^n,$ $|\psi(z)|\leq C$ for all $z\in K$ and $\phi(K)$ is also a
compact subset of $U^n,$ by the hypothesis,
$\left|\psi(z)f_{j_m}(\phi(z))\right|\leq
C\left|f_{j_m}(\phi(z))\right|$ converges to zero uniformly on
$K$. It follows from the arbitrary of $K$ that the limit function
$g$ is equal to $0.$ Since it is true for arbitrary subsequence of
$\{f_j\}$ we see that $W_{\psi,\phi}f_j\to 0$ in ${\cal B}^q.$

Conversely,  $\{g_j\}$ be any sequence in the ball ${\cal
K}_M=B_{{\cal B} ^p}(0,M)$ of the space ${\cal B} ^p(U^n).$ Since
$\|g_j\|_p\leq M<\infty,$ by Lemma 1, $\{g_j\}$ is uniformly
bounded on compact subsets of $U^n$ and hence normal by Montel's
theorem. Hence we may extract a subsequence $\{g_{j_m}\}$ which
converges uniformly on compact subsets of $U^n$ to some $g\in
H(U^n).$ It follows that $\frac{{\partial g}_{j_m}}{\partial
z_l}\to \frac{\partial g}{\partial z_l}$ uniformly on compacts,
for each $l\in\{1,...,n\},$ which implies $g\in {\cal B} ^p(U^n)$
and $\|g\|_p\leq M.$ Hence the sequence $\{g_{j_m}-g\}$ is such
that $\|g_{j_m}-g\|_p\leq 2M<\infty,$ and converges to $0$ on
compact subsets of $U^n,$ by the hypothesis and Lemma 2, we have
that $\psi g_{j_m}\circ\phi\to\psi g\circ\phi$ in ${\cal B}^q.$
Thus the set $W_{\psi,\phi}({\cal K}_M)$ is relatively compact,
finishing the proof.

\begin{Lemma} Let $0\leq p<1.$ If $f\in{\cal B} ^p(U^n)$, then
$$ |f(z)-f(w)|\leq \displaystyle\frac{2}{1-p}\|f\|_p\sum\limits^n_{k=1}|z_k-w_k|^{1-p},
\quad z,w\in U^n.$$\end{Lemma}

{\bf Proof}\hspace*{4mm} For any $z=(z_1,z_2,\cdots,z_n),$
$w=(w_1,w_2,\cdots,w_n)\in U^n,$ then $$tz+(1-t)w\in U^n$$ for
$t\in [0,1].$ Denote $F(t)=f\left(tz+(1-t)w\right),$ then
\begin{eqnarray}&&f(z)-f(w)=F(1)-F(0)=\int^1_0 F'(t)dt=\int^1_0\displaystyle\frac{d}{dt}\left(f(tz+(1-t)w\right)dt
\nonumber\\
&&=\sum\limits^n_{k=1}(z_k-w_k)\int^1_0\displaystyle\frac{\partial
f} {\partial
\zeta_k}\left(tz+(1-t)w\right)dt.\label{17}\end{eqnarray} $f\in
{\cal B} ^p(U^n),$ so for $k\in\{1,2,\cdots,n\}$
$$(1-|\zeta_k|^2)^p\left|\displaystyle\frac{\partial
f} {\partial \zeta_k}(\zeta)\right|\leq \|f\|_p$$
$$\left|\displaystyle\frac{\partial
f} {\partial \zeta_k}\right|\leq \|f\|_p(1-|\zeta_k|^2)^{-p}\leq
\|f\|_p(1-|\zeta_k|)^{-p}.$$ Let
$\zeta_k=tz_k+(1-t)w_k=w_k+t(z_k-w_k),$ note that if $|a|<1,|b|<1$
and $|a+b|<1$ then $|a|+|b|+|a+b|<2,$ so
$1-|\zeta_k|=1-|w_k+t(z_k-w_k)|\geq
\left|1-|w_k|-t|z_k-w_k|\right|,$
\begin{equation}\left|\displaystyle\frac{\partial f} {\partial
\zeta_k}(tz+(1-t)w)\right|\leq \|f\|_p(1-|\zeta_k|^2)^{-p}\leq
\|f\|_p\left|1-|w_k|-t|z_k-w_k|\right|^{-p}.\label{18}\end{equation}
It follows from (\ref{17}) and (\ref{18}) that
\begin{eqnarray}&&|f(z)-f(w)|\leq
\|f\|_p\sum\limits^n_{k=1}|z_k-w_k|\int^1_0\left|1-|w_k|-t|z_k-w_k|\right|^{-p}dt\nonumber\\
&&=\|f\|_p\sum\limits^n_{k=1}\int^{|z_k-w_k|}_0\left|1-|w_k|-t\right|^{-p}dt.\label{a}\end{eqnarray}
If $1-|w_k|\leq |z_k-w_k|,$ then
\begin{eqnarray}&&\int^{|z_k-w_k|}_0\left|1-|w_k|-t\right|^{-p}dt\nonumber\\
&&=\int^{1-|w_k|}_0(1-|w_k|-t)^{-p}dt+\int^{|z_k-w_k|}_{1-|w_k|}\left(t-(1-|w_k|)\right)^{-p}dt\nonumber\\
&&=\displaystyle\frac{1}{1-p}\left((1-|w_k|)^{1-p}+(|z_k-w_k|-(1-|w_k|))^{1-p}\right)\nonumber\\
&&\leq
\displaystyle\frac{2}{1-p}(|z_k-w_k|)^{1-p}.\label{b}\end{eqnarray}
If $1-|w_k|>|z_k-w_k|,$ let $a=1-|w_k|, b=|z_k-w_k|, a>b$ and note
that the fundamental inequality $a^{1-p}<(a-b)^{1-p}+b^{1-p},$
where $0\leq p<1,$ so
\begin{eqnarray}&&\int^{|z_k-w_k|}_0\left|1-|w_k|-t\right|^{-p}dt=\int^{|z_k-w_k|}_0\left(1-|w_k|-t\right)^{-p}dt\nonumber\\
&&=\displaystyle\frac{1}{1-p}\left((1-|w_k|)^{1-p}-(1-|w_k|-|z_k-w_k|)^{1-p}\right)\nonumber\\
&&<\displaystyle\frac{1}{1-p}|z_k-w_k|^{1-p}.\label{c}
\end{eqnarray}
Combining $(\ref{b}),(\ref{c})$ and $(\ref{a})$, this Lemma is
proved.

\begin{Lemma} Let $0\leq p<1.$ Every norm-bounded sequence in
${\cal B} ^p(U^n)$ has a subsequence that converges uniformly on
$\overline{U^n}$.\end{Lemma}

{\bf Proof}\hspace*{2mm} Let $\{f_j\}$ be a sequence in ${\cal
B}^p(U^n),$ and $\|f_j\|_p\leq C $ for $j=1,2,\cdots.$ It follows
from Lemma 3 that $|f_j(z)-f_j(w)|\leq
C\sum\limits^n_{k=1}|z_k-w_k|^{1-p}$ for $z,w\in U^n.$ Thus the
family $\{f_j: n=1,2,\cdots\}$ is equicontinuous. By Lemma 1,
$|f_j(z)|\leq C\|f_j\|_p\leq C, $ the family $\{f_j:
n=1,2,\cdots\}$ is bounded uniformly on $\overline{U^n}$. The
statement of the lemma now follows from the Arzela-Ascoli Theorem.

\begin{Lemma} Let $p\in [0,1)$ and $\{f_j\}$ be a norm bounded sequence in ${\cal B} ^p(U^n)$ which converges
to zero uniformly on compacts of $U^n.$ Then it converges to zero
uniformly on $\overline{U^n}.$
\end{Lemma}

{\bf Proof}\hspace*{4mm} Let $\varepsilon\in (0,1),$ $r\in
(1-\varepsilon,1),$ and $w\in rU^n,$ then from Lemma 3 we have
\begin{eqnarray*}&&|f_j(z)|\leq
|f_j(w)|+C\|f_j\|_p\sum\limits^n_{k=1}|z_k-w_k|^{1-p}\\
&&\leq |f_j(w)| +C\sum\limits^n_{k=1}|z_k-w_k|^{1-p},
\hspace*{4mm} z\in U^n.\end{eqnarray*} Since $\sup\limits_{w\in
rU^n}|f_j(w)|\to 0$ and since for each $z\in
\overline{U^n}\setminus{rU^n},$ $1-\varepsilon<r\leq |z_k|\leq 1,$
there exists $\eta$ such that
$$1-\varepsilon<\eta<r\leq |z_k|\leq 1.$$ Choosing
$w_z=(\eta z_1,\eta z_2,\cdots,\eta z_n),$ then $w_z\in {rU^n}$
and
$$|z_k-(w_z)_k|^{1-p}=|z_k-\eta z_k|^{1-p}=(1-\eta)^{1-p}|z_k|<\varepsilon^{1-p}\leq\varepsilon.$$
So $\lim\limits_{k\to\infty}\sup_{z\in \overline{U^n}}|f_j(z)|\leq
C\varepsilon,$ from which the result follows.

Combining Lemma 2 and Lemma 5, we can obtain the following Lemma
at once.

\begin{Lemma} Let $\phi$ be a holomorphic self-map of $U^n$ and $\psi(z)$ a holomorphic function
on $U^n.$ Let $0\leq p<1$ and $q\geq 0,$ then $W_{\psi,\phi}:{\cal
B}^p(U^n)\to {\cal B} ^q(U^n)$ is compact if and only if for any
bounded sequence $\{f_j\}$ in ${\cal B} ^p(U^n)$ which converges
to zero uniformly on $\overline{U^n}$ we have
$\|W_{\psi,\phi}f_j\|_q\to 0,$ as $j\to\infty.$
\end{Lemma}

\section{The Proof of Theorem 1}

First shows the sufficiency. If $f\in {\cal B}^p(U^n),$ a direct
calculation gives
\begin{eqnarray}&&\sum\limits^n_{k=1}\left|\displaystyle\frac{\partial
\left(W_{\psi,\phi}f(z)\right)} {\partial z_k}\right|
(1-|z_k|^2)^{q}\nonumber\\
&&=\sum\limits^n_{k=1}\left|\displaystyle\frac{\partial
\psi(z)}{\partial
z_k}f(\phi(z))+\psi(z)\sum\limits^n_{l=1}\displaystyle\frac{\partial
f}{\partial w_l}(\phi(z)) \displaystyle\frac{\partial
\phi_l}{\partial z_k}(z)\right|(1-|z_k|^2)^{q}\label{19}\\
&&\leq \sum\limits^n_{k=1}\left|\displaystyle\frac{\partial
\psi(z)}{\partial
z_k}\right||f(\phi(z))|(1-|z_k|^2)^{q}\nonumber\\
&& +|\psi(z)|\sum\limits^n_{k,
l=1}\left|\displaystyle\frac{\partial f}{\partial w_l}(\phi(z))
\displaystyle\frac{\partial \phi_l}{\partial
z_k}(z)\right|(1-|z_k|^2)^{q}.\label{20}\end{eqnarray} By Lemma 1,
if $p=1$, (\ref{19}) gives that
\begin{eqnarray}&&\sum\limits^n_{k=1}\left|\displaystyle\frac{\partial
\left(W_{\psi,\phi}f(z)\right)} {\partial z_k}\right|
(1-|z_k|^2)^{q}\nonumber\\
&&\leq
C\left(\sum\limits^n_{k,l=1}\left|\displaystyle\frac{\partial
\psi} {\partial z_k}(z)\right|\left(1-|z_k|^2\right)^{q}\ln
\displaystyle\frac{4}{1-|\phi_l(z)|^2}\right.\nonumber\\
&&\left.+|\psi(z)|\sum\limits^n_{k,l=1}\left|\displaystyle\frac{\partial
\phi_{l}} {\partial
z_k}(z)\right|\displaystyle\frac{\left(1-|z_k|^2\right)^{q}}{1-|\phi_l(z)|^2}
\right)\|f\|_{p}.\label{21}\end{eqnarray} Combining (\ref{1}),
(\ref{2}) and (\ref{21}), we know $W_{\psi,\phi}: {\cal
B}^p(U^n)\longrightarrow {\cal B}^q(U^n) $ is bounded.

If $0\leq p<1$, by Lemma 1, (\ref{20}) gives that
\begin{eqnarray}&&\sum\limits^n_{k=1}\left|\displaystyle\frac{\partial
\left(W_{\psi,\phi}f(z)\right)} {\partial z_k}\right|
(1-|z_k|^2)^{q}\nonumber\\
&&\leq
C\left(\sum\limits^n_{k,l=1}\left|\displaystyle\frac{\partial
\psi} {\partial
z_k}(z)\right|\left(1-|z_k|^2\right)^{q}\right.\nonumber\\
&&\left.+|\psi(z)|\sum\limits^n_{k,l=1}\left|\displaystyle\frac{\partial
\phi_{l}} {\partial
z_k}(z)\right|\displaystyle\frac{\left(1-|z_k|^2\right)^{q}}{\left(1-|\phi_l(z)|^2\right)^{p}}
\right)\|f\|_{p}.\label{22}\end{eqnarray} Combining (\ref{3}),
(\ref{4}), Lemma 1 and (\ref{22}), we know $W_{\psi,\phi}: {\cal
B}^p(U^n)\longrightarrow {\cal B}^q(U^n) $ is bounded.

If $p>1,$by Lemma 1, (\ref{20}) gives that
\begin{eqnarray}&&\sum\limits^n_{k=1}\left|\displaystyle\frac{\partial
\left(W_{\psi,\phi}f(z)\right)} {\partial z_k}\right|
(1-|z_k|^2)^{q}\nonumber\\
&&\leq
C\left(\sum\limits^n_{k,l=1}\left|\displaystyle\frac{\partial
\psi} {\partial
z_k}(z)\right|\displaystyle\frac{\left(1-|z_k|^2\right)^{q}}{\left(1-|\phi_l(z)|^2\right)^{p-1}}\right.\nonumber\\
&&\left.+|\psi(z)|\sum\limits^n_{k,l=1}\left|\displaystyle\frac{\partial
\phi_{l}} {\partial
z_k}(z)\right|\displaystyle\frac{\left(1-|z_k|^2\right)^{q}}{\left(1-|\phi_l(z)|^2\right)^{p}}
\right)\|f\|_{p}.\label{23}\end{eqnarray} Combining (\ref{5}),
(\ref{6}) and (\ref{23}), we know $W_{\psi,\phi}: {\cal
B}^p(U^n)\longrightarrow {\cal B}^q(U^n) $ is bounded.

To show the necessity, assume that $W_{\psi,\phi}: {\cal
B}^p(U^n)\longrightarrow {\cal B}^q(U^n) $ is bounded, with
\begin{equation}\|W_{\psi,\phi}f\|_{q}\leq C\|f\|_{p}\label{24}\end{equation}
for all $f\in {\cal B}^p(U^n).$ It is clear that $\psi\in{\cal
B}^q(U^n),$ that is, (\ref{3}) holds.

If $p=1$, For fixed $l (1\leq l\leq n),$ we will make use of a
family of test functions $\{f_{w}: w\in\mbox{\Bbb C}, |w|<1\}$ in
${\cal B}^p(U^n)$ defined as follows: let
$$f_w(z)=\displaystyle\frac{z_l}{1-\overline {w}z_l},$$ then $\displaystyle\frac{\partial
f_w(z) }{\partial z_l}=\displaystyle\frac{1}{(1-\overline
{w}z_l)^2},\hspace{4mm}\displaystyle\frac{\partial f_w(z)
}{\partial z_k}=0 \hspace*{4mm}(k\neq l)$,
\hspace*{4mm}$\|f_w\|_p\leq \displaystyle\frac{4}{1-|w|^2}.$

It follows from (\ref{19}) and (\ref{24}) that
$$\sum\limits^n_{k=1}\left|\displaystyle\frac{\partial
\psi(z)}{\partial z_k}\displaystyle\frac{\phi_l(z)}{1-\overline
{w}\phi_l(z)}+\psi(z)\displaystyle\frac{1}{(1-\overline
{w}\phi_l(z))^2}\displaystyle\frac{\partial \phi_l}{\partial
z_k}(z)\right|(1-|z_k|^2)^{q}\leq C
\displaystyle\frac{1}{1-|w|^2},$$let $w=\phi_l(z),$ then
\begin{equation}|\psi(z)|\sum\limits^n_{k=1}\left|\displaystyle\frac{\partial
\phi_{l}} {\partial
z_k}(z)\right|\displaystyle\frac{\left(1-|z_k|^2\right)^{q}}{1-|\phi_l(z)|^2}\leq
C+\|\psi\|_q\leq C.\label{25}\end{equation} So for any $z\in U^n$,
$$|\psi(z)|\sum\limits^n_{k,l=1}\left|\displaystyle\frac{\partial
\phi_{l}} {\partial
z_k}(z)\right|\displaystyle\frac{\left(1-|z_k|^2\right)^{q}}{1-|\phi_l(z)|^2}\leq
C,$$ that is, (\ref{2}) is valid.

If we set
$$g_w(z)=\ln\displaystyle\frac{4}{1-\overline
{w}z_l},$$ then $\displaystyle\frac{\partial g_w(z) }{\partial
z_l}=\displaystyle\frac{\overline w}{1-\overline
{w}z_l},\hspace{4mm}\displaystyle\frac{\partial g_w(z) }{\partial
z_k}=0 \hspace*{4mm}(k\neq l)$, \hspace*{4mm} $\|g_w\|_p\leq 2+\ln
4 .$

It follows from (\ref{19}) and (\ref{24}) that
$$\sum\limits^n_{k=1}\left|\displaystyle\frac{\partial
\psi(z)}{\partial z_k}\ln\displaystyle\frac{4}{1-\overline
{w}\phi_l(z)}+\psi(z)\displaystyle\frac{\overline
w}{1-\overline{w}\phi_l(z)}\displaystyle\frac{\partial
\phi_l}{\partial z_k}(z)\right|(1-|z_k|^2)^{q}\leq C,$$ let
$w=\phi_l(z)$, (\ref{1}) is valid.

If $0\leq p<1$, $\psi=W_{\psi,\phi}1\in{\cal B}^q(U^n),$ (3) is
valid. For fixed $l (1\leq l\leq n),$ we will make use of a family
of test functions $\{f_{w}: w\in\mbox{\Bbb C}-\{0\}, |w|<1\}$ in
${\cal B} ^p(U^n$ defined as follows: set
$$f_w(z)=\int^{z_l}_0\left(1-\displaystyle\frac{\overline w^2}{|w|^2}z_l^2
\right)^{-p}dz_l$$ with $f_w(0)=0.$ Since for $w\neq 0,$
$$\displaystyle\frac{\partial f_w}{\partial z_l}
=\left(1-\displaystyle\frac{\overline
w^2}{|w|^2}z_l^2\right)^{-p}, \hspace*{4mm}
\displaystyle\frac{\partial f_w}{\partial z_i}=0,
\hspace*{4mm}(i\neq l),$$ it is easy to show $\|f_w\|_{\alpha}=1.$
Lemma 1 gives $|f_w(z)|\leq 1+\displaystyle\frac{n}{1-p},$ so by
(\ref{3}),
\begin{equation}\sum\limits^n_{k=1}\left|\displaystyle\frac{\partial
\psi}{\partial z_k}(z)f_w(\phi(z)) (1-|z_k|^2)^q\right| \leq
C\sum\limits^n_{k=1} \left|\displaystyle\frac{\partial
\psi}{\partial z_k}(z)\right|\left(1-|z_k|^2 \right)^q\leq
C.\label{26}\end{equation} For $z\in U^n,$ it follows from
(\ref{19}) and (\ref{24}) that
$$\sum\limits^n_{k=1}\left|\displaystyle\frac{\partial
\psi}{\partial z_k}(z)f_w(\phi(z))
+\psi(z)\sum\limits^n_{l=1}\displaystyle\frac{\partial
f_w(\phi(z))}{\partial \phi_l} \displaystyle\frac{\partial
\phi_l}{\partial z_k}(z)\right|(1-|z_k|^2)^q\leq C,$$ combining
(\ref{26}), we get $$\sum\limits^n_{k=1}\left|\psi(z)
\displaystyle\frac{\partial f_w(\phi(z))}{\partial \phi_l}
\displaystyle\frac{\partial \phi_l}{\partial
z_k}(z)\right|(1-|z_k|^2)^q\leq C.$$ For an arbitrary $z\in U^n$
with $\phi_l(z)\neq 0$ and set $w=\phi_l(z)$ in the above
inequality to obtain
$$|\psi(z)|\sum\limits^n_{k=1}\left|\displaystyle\frac{\partial
\phi_{l}} {\partial z_k}(z)\right|
\displaystyle\frac{\left(1-|z_k|^2\right)^q}{\left(1-|\phi_l(z)|^2\right)^p}
\leq C.$$

If $\phi_l(z)=0,$ let $f(z)=z_l,$ it follows from (\ref{19}) and
(\ref{24}) that
$$|\psi(z)|\left|\displaystyle\frac{\partial
\phi_{l}} {\partial z_k}(z)\right|
\displaystyle\frac{\left(1-|z_k|^2\right)^q}{\left(1-|\phi_l(z)|^2\right)^p}
=\left|\displaystyle\frac{\partial\psi}{\partial z_k}(z)\phi_l(z)
+\psi(z)\displaystyle\frac{\partial\phi_l}{\partial
z_k}(z)\right|(1-|z_k|^2)^q \leq C.$$

So for all $z\in U^n,$
$$|\psi(z)|\sum\limits^n_{k,l=1}\left|\displaystyle\frac{\partial
\phi_{l}} {\partial z_k}(z)\right|
\displaystyle\frac{\left(1-|z_k|^2\right)^q}{\left(1-|\phi_l(z)|^2\right)^p}\leq
C,$$ that is, (\ref{4}) is valid.

If $p>1$, similarly, the test functions would be
$$f_w(z)=\displaystyle\frac{1}{\overline w}\left(\displaystyle\frac{1-|w|^2}{\left(1-\overline
{w}z_l\right)^p}-\displaystyle\frac{1}{\left(1-\overline
{w}z_l\right)^{p-1}}\right),\hspace*{4mm}(|w|<1, w\neq 0, 1\leq
l\leq n)$$ for proving (\ref{6}), and
$$g_w(z)=\left(\displaystyle\frac{p}{\left(1-\overline
{w}z_l\right)^{p-1}}-\displaystyle\frac{(p-1)(1-|w|^2)}{\left(1-\overline
{w}z_l\right)^p}\right),\hspace*{4mm}(|w|<1, 1\leq l\leq n)$$ for
proving (\ref{5}), we omit the details. Now the proof of Theorem 1
is completed.

\section{The Proof of Theorem 2}

To show the sufficiency.

If $p=1$. First assume conditions (\ref{7}) and (\ref{8}) hold, we
need to prove $W_{\psi,\phi}: {\cal B}^p(U^n)\longrightarrow {\cal
B}^q(U^n) $ is compact. According to Lemma 2, assume that
$\left\|f_{j}\right\|_{p}\leq C, j=1,2,\ldots,$ and $\{f_{j}\}$
converges to zero uniformly on compact subsets of $U^n,$ we need
only prove that $\left\|W_{\psi,\phi}f_{j}\right\|_{q}\to 0,$ as
$j\to\infty.$

For every $\varepsilon>0,$ (\ref{7}) and (\ref{8}) imply that
there exists a $r$, $0<r<1$, such that
\begin{equation}\sum\limits^n_{k,l=1}\left|\displaystyle\frac{\partial \psi}
{\partial z_k}(z)\right|\left(1-|z_k|^2\right)^{q}\ln
\displaystyle\frac{4}{1-|\phi_l(z)|^2} <\varepsilon
\label{27}\end{equation}
and\begin{equation}\sum\limits^n_{k,l=1}\left|\psi(z)\displaystyle\frac{\partial
\phi_l}{\partial z_k}(z)\right|
\displaystyle\frac{\left(1-|z_k|^2\right)^{q}}{1-|\phi_l(z)|^2}
<\varepsilon,\label{28}\end{equation} whenever $dist(\phi(z),
\partial U^n)<r.$
Since $f_j(\phi(0))$ converges to zero, by (\ref{27}) and
(\ref{28}), it follows from (\ref{21}) that for large enough $j$
\begin{equation}\sum\limits^n_{k=1}\left|\displaystyle\frac{\partial\left(W_{\psi,\phi}f_j\right)}
{\partial z_k}(z)\right| (1-|z_k|^2)^{q}\leq
C\varepsilon,\label{29}\end{equation} whenever $dist\left(\phi(z),
\partial U^n\right)< r.$ If we set $\varepsilon=1,$ by (\ref{27}) and
(\ref{28}),it is easy to show that (\ref{1}) and (\ref{2}) holds,
from Remark 1, we know $\psi\in{\cal B}^q(U^n).$

On the other hand, if we write $E=\left\{w\in U^n: dist(w,\partial
U^n)\geq r\right\},$ which is a closed subset of $U^n$, then
$f_j(w)$ and $\displaystyle\frac{\partial f_j}{\partial w_l}(w)\to
0$ uniformly on $E$. So by (\ref{19}) and the boundedness
condition (2) in Theorem 1, we obtain for $dist(\phi(z),\partial
U^n)\geq r,$
\begin{eqnarray}&&\sum\limits^n_{k=1}\left|\displaystyle\frac{\partial\left(W_{\psi,\phi}f_j\right)}
{\partial z_k}(z)\right| (1-|z_k|^2)^{q}\leq
C\left(|f_j(\phi(z))|+\sum\limits^n_{l=1}\left|\displaystyle\frac{\partial
f_j}{\partial w_l}(\phi(z))\right|\right)\nonumber\\
&& =
C\left(|f_j(w)|+\sum\limits^n_{l=1}\left|\displaystyle\frac{\partial
f_j}{\partial w_l}(w)\right|\right)\leq C\varepsilon.\label{30}
\end{eqnarray}

Since $\psi(0)f_j(\phi(0))\to 0$ for large enough $j,$ combining
(\ref{29}) and (\ref{30}) we know
$\left\|W_{\psi,\phi}f_{j}\right\|_{q}\to 0,$ as $j\to\infty.$

If $p>1,$ note that (\ref{23}), in a similar manner to the case
$p=1$, we can show the sufficiency, omit the details.

To show the necessity.

First we will prove the following Lemma: \begin{Lemma} Let $p\geq
1$ and $q\geq 0$. If $W_{\psi,\phi}: {\cal
B}^p(U^n)\longrightarrow {\cal B}^q(U^n)$ is compact, then
\begin{equation}
|\psi(z)|\sum\limits^n_{k,l=1}\left|\displaystyle\frac{\partial
\phi_{l}} {\partial
z_k}(z)\right|\displaystyle\frac{\left(1-|z_k|^2\right)^{q}}{\left(1-|\phi_l(z)|^2\right)^{p}}
=o(1)\hspace*{6mm} (\phi(z)\to\partial
U^n).\label{31}\end{equation}\end{Lemma}

{\bf Proof}\hspace{4mm} For any $\{z^{j}\}$ in $U^n$ with
$\phi(z^j)\to\partial U^n$ as $j\to\infty.$ Let
$A_j=\phi(z^{j})=(a_{1j},\ldots, a_{nj}),$ where
$a_{kj}=\phi_{k}(z^j), 1\leq k\leq n .$ Since
$\phi(z^{j})\to\partial U^n,$ as $j\to\infty$, there exists some
$s, (1\leq s\leq n),$ with $|a_{sj}|\to 1$ as $j\to\infty.$
Without loss of generality, we may assume that $s=1$.

Case 1: If for some $l (1\leq l\leq n),$ $|a_{lj}|\to 1$ as
$j\to\infty,$ then set
\begin{equation}f_{j}(z)=\displaystyle\frac{\left(1-|a_{lj}|^2\right)^2}{\left(1-z_l\overline{a_{lj}}\right)^{p+1}}-
\displaystyle\frac{1-|a_{lj}|^2}{\left(1-z_l\overline{a_{lj}}\right)^p},
\hspace*{4mm}j=1,2,\cdots.\label{32}
\end{equation}
$$\displaystyle\frac{\partial f_{j}}{\partial
z_{k}}=0,\hspace*{4mm} (k\neq l),$$ $$\displaystyle\frac{\partial
f_{j}}{\partial z_{l}}
=(p+1)\overline{a_{lj}}\displaystyle\frac{\left(1-|a_{lj}|^2\right)^2}{\left(1-z_l\overline{a_{lj}}\right)^{p+2}}-
p\overline{a_{lj}}\displaystyle\frac{1-|a_{lj}|^2}{\left(1-z_l\overline{a_{lj}}\right)^{p+1}},$$
so $\|f_{j}\|_{p}\leq C$ and $\{f_j(z)\}$ tends to zero uniformly
on compact subsets of $U^n.$ So by Lemma 2, we know
$\left\|W_{\psi,\phi}f_{j}\right\|_{q}\to 0$ as $j\to\infty.$ It
follows from (\ref{19}) that
\begin{eqnarray}&&\sum\limits^n_{k=1}\left(1-|z^j_k|^2\right)^{q}
\left|\displaystyle\frac{\partial \psi}{\partial
z_k}(z^j)f_j(\phi(z^j))+\psi(z^j)\displaystyle\frac{\partial
f_j}{\partial w_l}(\phi(z^j))\displaystyle\frac{\partial
\phi_l}{\partial z_k}(z^j)\right|\nonumber\\
&&\leq \left\|W_{\psi,\phi}f_{j}\right\|_{q}\to 0
.\label{33}\end{eqnarray} Note that $f_j(\phi(z^j))=0$ and
$\displaystyle\frac{\partial f_j}{\partial
z_l}(\phi(z^j))=\overline{a_{lj}}\displaystyle\frac{1}{\left(1-|a_{lj}|^2\right)^p}$,
(\ref{33}) gives that
\begin{eqnarray}&&|\psi(z^j)|\sum\limits^n_{k=1}
\displaystyle\frac{\left(1-|z^j_k|^2\right)^{q}}{\left(1-|a_{lj}|^2\right)^p}
\left|\displaystyle\frac{\partial
\phi_l}{\partial z_k}(z^j)\right|\nonumber\\
&&\leq
C\displaystyle\frac{1}{|a_{lj}|}\left\|W_{\psi,\phi}f_{j}\right\|_{q}
\to 0 \hspace*{4mm}\mbox{as $j\to\infty.$}\label{34}\end{eqnarray}

Case 2: If for some $l (2\leq l\leq n),$$|a_{lj}|\not\to 1$ as
$j\to\infty,$ then we assume $|a_{lj}|\leq\rho<1.$

If $p\geq 1$, The compactness of $W_{\psi,\phi}$ implies
$W_{\psi,\phi}$ is bounded, by Theorem 1 and Remark 1, (\ref{11})
is valid. Set
\begin{equation}f_{j}(z)=z_l\left(\displaystyle\frac{1-|a_{1j}|^2}
{1-z_1\overline{a_{1j}}}\right)^p,
\hspace*{4mm}j=1,2,\cdots.\label{35}
\end{equation}
It is easy to show that $\|f_j\|_{p}\leq C$, and $\{f_j\}$ tends
to zero uniformly on compact subsets of $U^n.$ It follows from
(\ref{19}) that
\begin{eqnarray*}&&|\psi(z^j)|\sum\limits^n_{k=1}
\left(1-|z^j_k|^2\right)^{q}\left|\displaystyle\frac{\partial
\phi_l}{\partial z_k}(z^j)\right|\nonumber\\
&&\leq \left\|W_{\psi,\phi}f_{j}\right\|_{q}
+C\sum\limits^n_{k=1}\left(1-|z^j_k|^2\right)^{q}\left|\displaystyle\frac{\partial
\psi }{\partial z_k}(z^j)\right|\nonumber\\
&&+C|\psi(z^j)|\sum\limits^n_{k=1}
\displaystyle\frac{\left(1-|z^j_k|^2\right)^{q}}{1-|a_{1j}|^2}
\left|\displaystyle\frac{\partial \phi_1}{\partial
z_k}(z^j)\right|
\nonumber\\
&&\leq \left\|W_{\psi,\phi}f_{j}\right\|_{q}
+C\sum\limits^n_{k=1}\left(1-|z^j_k|^2\right)^{q}\left|\displaystyle\frac{\partial
\psi }{\partial z_k}(z^j)\right|\nonumber\\
&&+C|\psi(z^j)|\sum\limits^n_{k=1}
\displaystyle\frac{\left(1-|z^j_k|^2\right)^{q}}{(1-|a_{1j}|^2)^p}
\left|\displaystyle\frac{\partial \phi_1}{\partial
z_k}(z^j)\right|\to 0.\label{36}\end{eqnarray*} So
\begin{eqnarray}&&|\psi(z^j)|\sum\limits^n_{k=1}
\displaystyle\frac{\left(1-|z^j_k|^2\right)^{q}}
{(1-|a_{lj}|^2)^p}\left|\displaystyle\frac{\partial
\phi_l}{\partial z_k}(z^j)\right|\nonumber\\
&&\leq
\displaystyle\frac{1}{(1-\rho^2)^p}|\psi(z^j)|\sum\limits^n_{k=1}
\left(1-|z^j_k|^2\right)^{q}\left|\displaystyle\frac{\partial
\phi_l}{\partial z_k}(z^j)\right|\to 0.\label{37}\end{eqnarray}

Now we return to prove the necessary of Theorem 2.

If $p=1$, by Lemma 7, (\ref{8}) is necessary, now we prove that
(\ref{7}) is also necessary.

In case 1,
set\begin{equation}f_j(z)=\left(\ln\displaystyle\frac{4}{1-|a_{lj}|^2}\right)^{-1}
\left(\ln\displaystyle\frac{4}{1-z_l\overline{a_{lj}}}\right)^2,\label{38}\end{equation}
It is easy to show that $\|f_j\|_{p}\leq C$, and $\{f_j\}$ tends
to zero uniformly on compact subsets of $U^n.$ It follows from
(\ref{19}) that
$$\sum\limits^n_{k=1}\left(1-|z^j_k|^2\right)^{q}
\left|\displaystyle\frac{\partial \psi}{\partial
z_k}(z^j)f_j(\phi(z^j))+\psi(z^j)\displaystyle\frac{\partial
f_j}{\partial w_l}(\phi(z^j))\displaystyle\frac{\partial
\phi_l}{\partial z_k}(z^j)\right|\leq
\left\|W_{\psi,\phi}f_{j}\right\|_{q}\to 0 .$$ By a direct
calculation and by (\ref{31}), we obtain
\begin{eqnarray}&&\sum\limits^n_{k=1}\left|\displaystyle\frac{\partial \psi}
{\partial z_k}(z)\right|\left(1-|z_k|^2\right)^{q}\ln
\displaystyle\frac{4}{1-|\phi_l(z)|^2}\nonumber\\
&&\leq
\left\|W_{\psi,\phi}f_{j}\right\|_{q}+2|\psi(z^j)|\sum\limits^n_{k=1}
\displaystyle\frac{\left(1-|z^j_k|^2\right)^{q}}{\left(1-|a_{lj}|^2\right)^{p}}
\left|\displaystyle\frac{\partial \phi_l}{\partial
z_k}(z^j)\right|\to 0.\label{39}\end{eqnarray}

In case 2, by remark 1,
$$\sum\limits^n_{k=1}\left|\displaystyle\frac{\partial \psi}
{\partial z_k}(z^j)\right|\left(1-|z^j_k|^2\right)^{q}\ln
\displaystyle\frac{4}{1-|\phi_l(z^j)|^2}\leq
\ln\displaystyle\frac{4}{1-\rho^2}\sum\limits^n_{k=1}\left|\displaystyle\frac{\partial
\psi} {\partial z_k}(z^j)\right|\left(1-|z^j_k|^2\right)^{q}\to
0.$$ Combining case 1 and case 2, we know (\ref{7}) is necessary.

If $p>1.$ by Lemma 7, we know (\ref{10}) is necessary, now we
prove (\ref{9}) is also necessary.

In case 1, set
$$f_{j}(z)=(p+1)\displaystyle\frac{1-|a_{lj}|^2}{\left(1-z_l\overline{a_{lj}}\right)^p}-
p\displaystyle\frac{\left(1-|a_{lj}|^2\right)^2}{\left(1-z_l\overline{a_{lj}}\right)^{p+1}},
\hspace*{4mm}j=1,2,\cdots.$$ It is easy to show that
$\|f_j\|_{p}\leq C$, and $\{f_j\}$ tends to zero uniformly on
compact subsets of $U^n.$ Note that
$f_j(\phi(z^j))=\displaystyle\frac{1}{\left(1-|a_{lj}|^2\right)^{p-1}}$
and $\displaystyle\frac{\partial f_j}{\partial z_l}(\phi(z^j))=0,$
it follows from (\ref{19}) that
\begin{eqnarray}&&\sum\limits^n_{k=1}\left(1-|z^j_k|^2\right)^{q}
\left|\displaystyle\frac{\partial \psi}{\partial
z_k}(z^j)f_j(\phi(z^j))+\psi(z^j)\displaystyle\frac{\partial
f_j}{\partial w_l}(\phi(z^j))\displaystyle\frac{\partial
\phi_l}{\partial
z_k}(z^j)\right|\nonumber\\
&&=\sum\limits^n_{k=1}\left|\displaystyle\frac{\partial
\psi}{\partial
z_k}(z^j)\right|\displaystyle\frac{\left(1-|z^j_k|^2\right)^{q}}{(1-|a_{lj}|^2)^{p-1}}\leq
\left\|W_{\psi,\phi}f_{j}\right\|_{q}\to
0.\label{40}\end{eqnarray}

In case 2, by remark 1,
\begin{eqnarray*}&&\sum\limits^n_{k=1}\left|\displaystyle\frac{\partial \psi}
{\partial
z_k}(z^j)\right|\displaystyle\frac{\left(1-|z^j_k|^2\right)^{q}}{\left(1-|\phi_l(z^j)|^2\right)^{p-1}}\\
&&\leq
\displaystyle\frac{1}{\left(1-\rho^2\right)^{p-1}}\sum\limits^n_{k=1}\left|\displaystyle\frac{\partial
\psi} {\partial z_k}(z^j)\right|\left(1-|z^j_k|^2\right)^{q} \to
0.\end{eqnarray*} Combining cas1 and case 2, we know (\ref{9}) is
necessary.

Now the proof of Theorem 2 is completed.

\section{Proof of Theorem 3}

Proof of Theorem 3. First assume $W_{\psi,\phi}: {\cal
B}^p(U^n)\longrightarrow {\cal B}^q(U^n)$ is bounded and
(\ref{12}) holds, we prove that $W_{\psi,\phi}:{\cal B} ^p(U^n)\to
{\cal B} ^q(U^n)$ is compact. By Lemma 6, assume that
$\left\|f_{j}\right\|_p\leq C,\, j=1,2,...,$ and $\{f_{j}\}$
converges to zero uniformly on $\overline{U^n},$ we need to prove
that $\left\|W_{\psi,\phi}f_{j}\right\|_q\to 0,$ as $j\to\infty.$
Note that the boundedness of $W_{\psi,\phi}$ and Theorem 1,
(\ref{1}) and (\ref{2}) hold. So
\begin{equation}|\psi(z)|\leq C\|\psi\|_q\leq C,\hspace*{4mm} |\psi(z)|\sum\limits^n_{k=1}\left|\frac{\partial
\phi_l}{\partial
z_k}(z)\right|\frac{(1-|z_k|^2)^q}{(1-|\phi_l(z)|^2)^p}\leq C
\label{41}\end{equation} for $l\in\{1,2,\cdots,n\}.$

By the assumption as $j\to\infty$
\begin{equation}\sup\limits_{z\in\overline{U^n}}|f_j(z)|\to 0\label{42}\end{equation}
\begin{eqnarray}&&\sum\limits^n_{k=1}\left|\displaystyle\frac{\partial
\left(W_{\psi,\phi}f_j(z)\right)} {\partial z_k}\right|
(1-|z_k|^2)^{q}\nonumber\\
&&=\sum\limits^n_{k=1}\left|\displaystyle\frac{\partial
\psi(z)}{\partial
z_k}f_j(\phi(z))+\psi(z)\sum\limits^n_{l=1}\displaystyle\frac{\partial
f_j}{\partial w_l}(\phi(z)) \displaystyle\frac{\partial
\phi_l}{\partial z_k}(z)\right|(1-|z_k|^2)^{q}\nonumber\\
&&\leq
|f_j(\phi(z))|\sum\limits^n_{k=1}\left|\displaystyle\frac{\partial
\psi} {\partial
z_k}(z)\right|\left(1-|z_k|^2\right)^{q}\nonumber\\
&&+\sum\limits^n_{k=1}\sum\limits^n_{l=1}\left|\frac{\partial
f_j}{\partial w_l}(\phi(z))\right||\psi(z)|\left|\frac{\partial
\phi_l} {\partial
z_k}(z)\right|\displaystyle\frac{(1-|z_k|^2)^q}{(1-|\phi_l(z)|^2)^p}\nonumber\\
&&\leq
C\left|f_j(\phi(z))\right|+\sum\limits^n_{k=1}\sum\limits^n_{k=1}\left|\frac{\partial
f_j}{\partial w_l}(\phi(z))\right||\psi(z)|\left|\frac{\partial
\phi_l} {\partial
z_k}(z)\right|\displaystyle\frac{(1-|z_k|^2)^q}{(1-|\phi_l(z)|^2)^p}\nonumber\\
&&\leq
C\varepsilon+\sum\limits^n_{k=1}\sum\limits^n_{l=1}\left|\frac{\partial
f_j}{\partial w_l}(\phi(z))\right||\psi(z)|\left|\frac{\partial
\phi_l} {\partial
z_k}(z)\right|\displaystyle\frac{(1-|z_k|^2)^q}{(1-|\phi_l(z)|^2)^p}.\label{43}\end{eqnarray}

For every $\varepsilon>0$, $l\in\{1,...,n\},$ (\ref{12}) implies
that there exists an $r$, $0<r<1$, such that
\begin{equation}|\psi(z)|\sum\limits^n_{k=1}\left|\frac{\partial
\phi_l}{\partial
z_k}(z)\right|\frac{(1-|z_k|^2)^q}{(1-|\phi_l(z)|^2)^p}<\varepsilon,\label{44}\end{equation}
whenever $|\phi_l(z)|>r.$

Note that $$U^n= \left\{z\in U^n:
|\phi_l(z)|>r\right\}\bigcup\left\{z\in U^n: |\phi_l(z)|\leq
r\right\}.$$

For each $l\in\{1,2,\cdots,n\}.$ If $|\phi_l(z)|>r,$ from
(\ref{44}) and $\|f_j\|_p\leq C$, we have
\begin{equation}\left|\frac{\partial f_j}{\partial
w_l}(\phi(z))\right||\psi(z)|\left|\frac{\partial \phi_l}
{\partial
z_k}(z)\right|\displaystyle\frac{(1-|z_k|^2)^q}{(1-|\phi_l(z)|^2)^p}<C\varepsilon.\label{45}\end{equation}

If $|\phi_l(z)|\leq r,$ by the Cauchy's estimate applied to the
function
$$g(w_l)=f(w_1,...,w_{l-1},w_l,w_{l+1}...,w_n),$$
of one variable, if $|w_l|\leq r,$ we have
$$\left|\frac{\partial f_j}{\partial w_l}(w)\right|\leq C\sup_{|w_l|\leq \frac{1+r}2}|f_j(w)|\leq C\, \sup_{z\in \overline{U^n}}|f_j(z)|,$$
for some $C>0$ independent of $f.$ From this, (\ref{41}) and
(\ref{42}) we have
\begin{eqnarray}&&\left|\frac{\partial
f_j}{\partial w_l}(\phi(z))\right||\psi(z)|\left|\frac{\partial
\phi_l} {\partial
z_k}(z)\right|\displaystyle\frac{(1-|z_k|^2)^q}{(1-|\phi_l(z)|^2)^p}\nonumber\\
&&\leq C\sup_{z\in \overline{U^n}}|f_j(z)|\to 0
\label{46}\end{eqnarray} as $j\to\infty.$ So

From (\ref{43}),(\ref{45}) and (\ref{46}), and since
$\lim\limits_{j\to\infty}f_j(\phi(0))=0,$ we obtain
$\left\|C_{\phi}f_{j}\right\|_q\to 0,$ as $j\to\infty,$ from which
the result follows.

Now suppose that $W_{\psi,\phi}: {\cal B} ^p(U^n)\to {\cal
B}^q(\mbox{\Bbb C}^n)$ is compact. Then $W_{\psi,\phi}$ is
bounded, and by Theorem 1 we know that (\ref{4}) holds. Now we
prove that condition (\ref{12}) holds. Let $\{z^{j}\}$ be a
sequence in $U^n$ and $w^j=\phi(z^{j})=(w^j_1,..., w^j_n).$ For
each $l\in\{1,2,\cdots,n\}$ with
$\lim\limits_{j\to\infty}|w^j_l|\to 1.$

Let
$$f_{j}(z)=\frac{1-|w^j_l|^2}{(1-z_l\overline{w^j_l})^p},\quad j=1,2,...\;,$$
then$$\displaystyle\frac{\partial f_j}{\partial
z_l}(z)=p\overline{w^j_l}\displaystyle\frac{1-|w^j_l|^2}{(1-z_l\overline{w^j_l|^2})^{p+1}}.$$
As in the proof of Theorem 2 we can prove that $\|f_{j}\|_p\leq
1+2^{p+1}p,$ for all $j\in {\bf N},$ i.e. $\{f_j\}$ is bounded on
${\cal B} ^p(U^{n}).$ It is easy to see that $f_j$ tends to zero
uniformly on compact subsets of $\mbox{\Bbb C}^n.$ So by Lemma 2,
we know $\left\|W_{\phi}f_{j}\right\|_q\to 0$ as $j\to\infty$ and
notice that
$$f_j(\phi(z^j))=(1-|w^j_l|^2)^{1-p},\hspace*{4mm}\displaystyle\frac{\partial
f_j}{\partial
w_l}(\phi(z^j))=p\overline{w^j_l}\displaystyle\frac{1}{(1-|w^j_l|^2)^p},$$
consequently, by (\ref{19})
\begin{eqnarray}&&\sum\limits^n_{k=1}\left(1-|z^j_k|^2\right)^{q}
\left|\displaystyle\frac{\partial \psi}{\partial
z_k}(z^j)f_j(\phi(z^j))+\psi(z^j)\displaystyle\frac{\partial
f_j}{\partial w_l}(\phi(z^j))\displaystyle\frac{\partial
\phi_l}{\partial z_k}(z^j)\right|\nonumber\\
&&\leq
\left\|W_{\psi,\phi}f_{j}\right\|_{q}.\label{47}\end{eqnarray} So
$$\sum\limits^n_{k=1}p|w^j_l||\psi(z^j)\left|\displaystyle\frac{\partial
\phi_l}{\partial
z_k}(z^j)\right|\displaystyle\frac{\left(1-|z^j_k|^2\right)^{q}}{(1-|w^j_l|^2)^p}\leq
\|\psi\|_q\left(1-|w^j_l|^2\right)^{1-p} +
\left\|W_{\psi,\phi}f_{j}\right\|_{q}.$$ Hence
$$|\psi(z^j)|\sum\limits^n_{k=1}\frac{(1-|z^j_k|^2)^q}{(1-|w^j_l|^2)^p}
\left|\frac{\partial \phi_l}{\partial z_k}(z^j)\right|\leq
\frac{1}{p|w^j_l|}\left(C(1-|w^j_l|^2)^{1-p} +
\left\|W_{\psi,\phi}f_{j}\right\|_{q}\right)\to 0$$ as
$j\to\infty.$ Now the proof of Theorem 3 is finished.

\end{document}